\documentclass[12pt,a4paper]{amsart}
\usepackage{amscd,amsfonts,amssymb,verbatim,mathrsfs}
\usepackage{amsthm}
\usepackage{amsmath}
\usepackage[dvips]{graphicx}
\usepackage{graphics,hyperref}
\usepackage{enumitem}
\usepackage{psfrag}
\usepackage{latexsym}

\def\frk{\frak}               

\def\Phi{{\frk n}}
\def\Phi{{\frk N}}
\def\opn#1#2{\def#1{\operatorname{#2}}} 
\opn\chara{char} \opn\length{\ell} \opn\pd{pd} \opn\rk{rk}
\opn\projdim{proj\,dim} \opn\injdim{inj\,dim} \opn\rank{rank}
\opn\depth{depth} \opn\grade{grade} \opn\height{height}
\opn\embdim{emb\,dim} \opn\codim{codim} \opn\sgn{sgn}

\opn\Tr{Tr} \opn\bigrank{big\,rank}
\opn\superheight{superheight}\opn\lcm{lcm}
\opn\trdeg{tr\,deg}
\opn\reg{reg} \opn\lreg{lreg} \opn\ini{in} \opn\lpd{lpd}
\opn\size{size}\opn\bigsize{bigsize}
\opn\cosize{cosize}\opn\bigcosize{bigcosize}
\opn\sdepth{sdepth}\opn\sreg{sreg}
\opn\link{link}\opn\fdepth{fdepth}
\opn\div{div} \opn\Div{Div} \opn\cl{cl} \opn\Cl{Cl} \opn\Cor{Cor}
\opn\Spec{Spec} \opn\Supp{Supp} \opn\supp{supp} \opn\Sing{Sing}
\opn\Ass{Ass} \opn\Min{Min}\opn\Mon{Mon} \opn\dstab{dstab} \opn\astab{astab}
\opn\Ann{Ann} \opn\Rad{Rad} \opn\Soc{Soc} \opn\Gr{Gr}
\opn\Im{Im} \opn\Ker{Ker} \opn\Coker{Coker} \opn\Am{Am}
\opn\Hom{Hom} \opn\Tor{Tor} \opn\Ext{Ext} \opn\End{End}
\opn\Aut{Aut} \opn\id{id} \opn\span{span}

\opn\nat{nat}
\opn\pff{pf}
\opn\Pf{Pf} \opn\GL{GL} \opn\SL{SL} \opn\mod{mod} \opn\ord{ord}
\opn\Gin{Gin} \opn\Hilb{Hilb}\opn\sort{sort} \opn\Gale{Gale}
\opn\aff{aff} \opn\conv{conv} \opn\relint{relint} \opn\st{st}   \opn\cone{cone}
\opn\lk{lk} \opn\cn{cn} \opn\core{core} \opn\vol{vol}
\opn\link{link} \opn\star{star}\opn\lex{lex} \opn\Gr{Gr}
\opn\gr{gr}

\def\pot#1#2{#1[\kern-0.28ex[#2]\kern-0.28ex]}

\opn\dirlim{\underrightarrow{\lim}}
\opn\inivlim{\underleftarrow{\lim}}
\def\Implies{\ifmmode\Longrightarrow \else
        \unskip${}\Longrightarrow{}$\ignorespaces\fi}
\def\implies{\ifmmode\Rightarrow \else
        \unskip${}\Rightarrow{}$\ignorespaces\fi}
\def\iff{\ifmmode\Longleftrightarrow \else
        \unskip${}\Longleftrightarrow{}$\ignorespaces\fi}

\let\:=\colon
\newtheorem{Theorem}{Theorem}[section]

\newtheorem{Corollary}[Theorem]{Corollary}
\newtheorem{Proposition}[Theorem]{Proposition}
\newtheorem{Remark}[Theorem]{Remark}

\newtheorem{Example}[Theorem]{Example}

\newtheorem{Definition}[Theorem]{Definition}

\let\epsilon\varepsilon
\let\kappa=\varkappa
\makeatletter

\def\qed{\ifhmode\textqed\fi
      \ifmmode\ifinner\quad\qedsymbol\else\dispqed\fi\fi}
\def\textqed{\unskip\nobreak\penalty50
       \hskip2em\hbox{}\nobreak\hfil\qedsymbol
       \parfillskip=0pt \finalhyphendemerits=0}
\def\dispqed{\rlap{\qquad\qedsymbol}}

\opn\dis{dis}
\def\pnt{{\raise0.5mm\hbox{\large\bf.}}}

\opn\Lex{Lex}

\begin{document}


\title{Radical splittings of toric ideals}

\author[1]{Anargyros Katsabekis}
\thanks{Corresponding author: Anargyros Katsabekis}

\author[2]{Apostolos Thoma}

\address{Anargyros Katsabekis, Department of Mathematics, University of Ioannina, Ioannina 45110, Greece}
\email{katsampekis@uoi.gr}

\address{Apostolos Thoma, Department of Mathematics, University of Ioannina, Ioannina 45110, Greece}
\email{athoma@uoi.gr}


\begin{abstract} Let $I_A \subset K[x_1,\ldots,x_n]$ be a toric ideal. In this paper, we provide a necessary and sufficient condition for the toric variety $V(I_A)$, over an algebraically closed field, to be expressed as the set-theoretic intersection of other toric varieties. We also introduce the radical splitting number of $I_A$, denoted by ${\rm Split}_{\rm rad}(I_A)$, and compute its exact value in several cases, with particular emphasis on toric ideals arising from graphs. In particular, we show that ${\rm Split}_{\rm rad}(I_A)=3$ for toric ideals of complete bipartite graphs. Additionally, we prove that ${\rm Split}_{\rm rad}(I_A)$ coincides with the binomial arithmetical rank of $I_A$ when the height of $I_A$ is equal to 2. 

\end{abstract}
\maketitle

\section{Introduction}

A natural question in algebraic geometry is to determine when an affine algebraic variety $V \subset K^n$ can be expressed as the intersection of two distinct affine algebraic varieties, $V_1$ and $V_{2}$, over an algebraically closed field $K$; that is, $V=V_{1} \cap V_{2}$, with $V \subsetneqq V_i$ for $1 \leq i \leq 2$. This geometric question translates algebraically via the Hilbert Nullstellensatz into an ideal-theoretic question: $$I(V)={\rm rad}(I(V_1)+I(V_2)),$$ where $I(W)$ denotes the ideal of an affine algebraic variety $W$ and ${\rm rad}(J)$ is the radical of an ideal $J$, see \cite[Chapter 4]{CLO}. This paper addresses the latter problem in the case where all varieties involved are toric varieties. In other words, we study when a given toric ideal $I_{A} \subset K[x_1,\ldots,x_n]$ can be written as $I_{A}={\rm rad}(I_{A_1}+I_{A_2})$ and $I_{A_i} \subsetneqq I_A$ for $1 \leq i \leq 2$.
Recent advances in understanding toric splittings \cite{FHKT, GS, KT2, KT3}, along with several earlier results on the binomial generation of toric ideals up to radical \cite{Katsa2, KMT}, allow us to provide a necessary and sufficient condition for the equality $I_{A}={\rm rad}(I_{A_1}+I_{A_2})$.

Toric ideals arising from combinatorial objects, and in particular from graphs, provide a rich framework in which such decompositions can be studied explicitly. In this paper, we place special emphasis on toric ideals of bipartite graphs and complete bipartite graphs, where the interaction between algebraic and combinatorial structure leads to precise and
computable invariants.

Let $A=\{{\bf a}_{1},\ldots,{\bf a}_{n}\}$ be a vector configuration in $\mathbb{Z}^{m}$, ${\rm ker}_{\mathbb{Q}}(A)=\{{\bf u}=(u_{1},\ldots,u_{n}) \in \mathbb{Q}^{n} \mid u_{1}{\bf a}_{1}+\cdots+u_{n}{\bf a}_{n}={\bf 0}\}$, and ${\rm ker}_{\mathbb{Z}}(A)={\rm ker}_{\mathbb{Q}}(A) \cap \mathbb{Z}^{n}$. Throughout this paper, we will assume that the affine semigroup $\mathbb{N}A=\{l_{1}{\bf a}_{1}+\cdots+l_{n}{\bf a}_{n} \mid l_{i} \in \mathbb{N}\}$ is pointed; namely ${\rm ker}_{\mathbb{Z}}(A) \cap \mathbb{N}^{n}=\{{\bf 0}\}$. Consider the polynomial ring $K[x_1,\ldots,x_n]$ over a field $K$. The {\em toric ideal} $I_A$ is the kernel of the $K$-algebra homomorphism $\phi: K[x_{1},\ldots,x_{n}] \longrightarrow K[t_{1},\ldots,t_{m},t_{1}^{-1},\ldots,t_{m}^{-1}]$ given by $\phi(x_{i})={\bf t}^{{\bf a}_{i}}$ for all $i=1,\ldots,n$ and it is generated by all the binomials $B({\bf u}):={\bf x}^{{\bf u}^{+}}-{\bf x}^{{\bf u}^{-}}$ such that ${\bf u}^{+}-{\bf u}^{-} \in {\rm ker}_{\mathbb{Z}}(A)$, where ${\bf u}^{+} \in \mathbb{N}^{n}$ and ${\bf u}^{-} \in \mathbb{N}^{n}$ denote the positive and negative part of ${\bf u}={\bf u}^{+}-{\bf u}^{-}$, respectively (see, for instance, \cite[Lemma 4.1]{St}). By \cite[Lemma 4.2]{St}, the height ${\rm ht}(I_A)$ of $I_A$ equals the dimension of the vector space ${\rm ker}_{\mathbb{Q}}(A)$. The set $V(I_A)=\{(u_{1},\ldots,u_{n}) \in K^n \mid F(u_1,\ldots,u_{n})=0, \forall F \in I_A\}$
of zeroes of $I_A$ is called {\em affine toric variety}, which is not necessarily normal.

The toric ideal $I_A$ is {\em splittable} if there exist toric ideals $I_{A_1}, I_{A_2}$ such that $I_A=I_{A_1}+I_{A_2}$ and $I_{A_i}\not =I_{A}$ for all $1 \leq i \leq 2$. The ideal $I_A$ is called {\em radical splittable} if and only if there exist toric ideals $I_{A_1}, I_{A_2}$ such that $I_{A}={\rm rad}(I_{A_1}+I_{A_2})$ and $I_{A_i} \neq I_A$, for $1 \leq i \leq 2$. 

\begin{Definition}
Let $I_A$ be a toric ideal.
\begin{enumerate}
    \item The {\em splitting number} of $I_A$, denoted by ${\rm Split}(I_A)$, is the smallest integer $s$ such that $I_{A} = I_{A_1} + \cdots + I_{A_s}$ and $I_{A_i} \neq I_A$ for every $1 \leq i \leq s$. In this case, $I_{A} = I_{A_1} + \cdots + I_{A_s}$ is a splitting of $I_A$.
\item The {\em radical splitting number} of $I_A$, denoted by ${\rm Split}_{\mathrm{rad}}(I_A)$, is the smallest integer $r$ such that $I_A = {\rm rad}(I_{A_1} + \cdots + I_{A_r})$ and $I_{A_i} \neq I_A$ for every $1 \leq i \leq r$. In this case, $I_{A} = {\rm rad}(I_{A_1} + \cdots + I_{A_r})$ is a radical splitting of $I_A$.
\end{enumerate}
\end{Definition}

\begin{Remark} {\rm  (1) Since $I_A$ is minimally generated by binomials, we get that ${\rm Split}(I_A) \leq \mu(I_A)$, where $\mu(I_A)$ is the minimal number of generators of $I_A$.\\ (2) The invariant ${\rm Split}(I_A)$ is an upper bound for ${\rm Split}_{\mathrm{rad}}(I_A)$.}
\end{Remark}

The {\em binomial arithmetical rank} of $I_A$, denoted by ${\rm bar}(I_A)$, is the smallest integer $t$ for which there exist binomials $B_{1},\ldots,B_{t}$ in $I_A$ such that $I_{A}={\rm rad}(B_{1},\ldots,B_{t})$. From the generalized Krull’s principal ideal theorem, we deduce that ${\rm ht}(I_A)$ is a lower bound for ${\rm bar}(I_A)$. 
From the definitions, we obtain the following inequalities for a non-principal toric ideal $I_A$: $$2 \leq {\rm Split}_{\mathrm{rad}}(I_A) \leq {\rm bar}(I_A) \leq \mu(I_A).$$ For a non-principal toric ideal $I_A$ we have that ${\rm Split}_{\mathrm{rad}}(I_A)=2$ if and only if $I_A$ is radical splittable.

In Section 2, we provide a necessary and sufficient condition for a toric ideal $I_A$ to be radical splittable (see Corollary \ref{main}). Also, we show that if the height of $I_A$ is $r \geq 3$ and ${\rm bar}(I_A) \leq 2r - 2$, then $I_A$ is radical splittable (see Theorem \ref{BasicMinimal}). Moreover, if $I_A$ is a toric ideal of height $r \geq 3$ such that $V(I_A)$ is a simplicial toric variety with full parametrization, then $I_A$ is radical splittable (see Proposition \ref{simplicial}). 

In Section 3, we prove that for a bipartite graph $G$, ${\rm Split}(I_G) = {\rm Split}_{\rm rad}(I_G)$ holds, where $I_G$ is the toric ideal of the graph $G$ (see Theorem \ref{Indispensable-Split}). Furthermore, we show that ${\rm Split}_{\rm rad}(I_G) = 3$ for a complete bipartite graph $G$ (see Theorem \ref{BipartiteSplitting}).

In Section 4, we study the special case in which $I_A$ is a toric ideal of height 2. We show that ${\rm Split}(I_A)=\mu(I_A)$ and ${\rm Split}_{\rm rad}(I_A)={\rm bar}(I_A)$ (see Theorem \ref{Radicalheight2}). We also show (see Corollary \ref{Simplicial-2}) that if $V(I_A)$ is a simplicial toric variety, then $2 \leq {\rm Split}_{\rm rad}(I_A) \leq 3$, and both lower and upper bounds are attained; see Example \ref{Simplicial-Example}. Furthermore, we explicitly compute the splitting and the radical splitting number of the toric ideal $I_{\Lambda(A)}$ (see Proposition \ref{Lawrence}), where $\Lambda(A)$ is the Lawrence lifting of the vector configuration $A$ associated with a symmetric monomial curve in $\mathbb{P}^{3}$. Proposition \ref{Lawrence} shows that the splitting number can be arbitrarily large for some configurations, while the radical splitting number remains constant for them. Finally, we show that both ${\rm Split}(I_A)$ and ${\rm Split}_{\rm rad}(I_A)$ can be arbitrarily large for certain other configurations; see Proposition \ref{cyclic}.

\section{Radical Splitting Criterion for $I_A$}
\label{section}

Let $A=\{{\bf a}_{1},\ldots,{\bf a}_{n}\}$ be a vector configuration in $\mathbb{Z}^{m}$. Given a set $C=\{{\bf c}_{1},\ldots,{\bf c}_{k}\}\subset \mathbb{Z}^n$, we shall denote by ${\rm span}_{\mathbb{Q}}(C)$ the $\mathbb{Q}$-vector space generated by the vectors of $C$, namely, $${\rm span}_{\mathbb{Q}}(C)=\{\lambda_{1}{\bf c}_{1}+\cdots+\lambda_{k}{\bf c}_{k} \mid \lambda_{i} \in \mathbb{Q}\}.$$ 

A set $S$ is a {\em minimal system of binomial generators} of the toric ideal $I_A$ {\em up to radical} if there exists no $S' \subsetneqq S$ such that $S'$ generates $I_A$ up to radical.

\begin{Theorem} \label{main2} There are toric ideals $I_{A_i}$ with $I_{A_i} \neq I_{A}$ for every $1 \leq i \leq s$, such that  $I_{A}=I_{A_1}+\cdots+I_{A_s}$ (respectively, $I_{A}={\rm rad}(I_{A_1}+\cdots+I_{A_s})$) if and only if there exists a minimal system of binomial generators $\{B({\bf u}) \mid {\bf u}\in C\subset {\rm ker}_{\mathbb{Z}}(A)\}$ of the toric ideal $I_A$ (respectively, of $I_A$ up to radical), and sets $C_{i}$, $1 \leq i \leq s$, such that $C=\cup_{i=1}^{s}C_{i}$ and ${\rm span}_{\mathbb{Q}}(C_i)\subsetneqq {\rm ker}_{\mathbb{Q}}(A)$ for every $1 \leq i \leq s$.  
\end{Theorem}

{\em \noindent Proof.} Suppose that $I_{A}=I_{A_1}+\cdots+I_{A_s}$ (respectively, $I_{A}={\rm rad}(I_{A_1}+\cdots+I_{A_s})$) and $I_{A_i} \neq I_A$, for every $1 \leq i \leq s$. Thus, ${\rm ker}_{\mathbb{Z}}(A_i)\subsetneqq {\rm ker}_{\mathbb{Z}}(A)$; therefore, ${\rm ker}_{\mathbb{Q}}(A_i)\subsetneqq {\rm ker}_{\mathbb{Q}}(A)$ for every $1 \leq i \leq s$.
Let $T_{i}$, $1 \leq i \leq s$, be a minimal system of binomial generators of $I_{A_i}$ (respectively, of $I_{A_i}$ up to radical). Then, since $I_{A}=I_{A_1}+\cdots+I_{A_s}$ (respectively, $I_{A}={\rm rad}(I_{A_1}+\cdots+I_{A_s})$), we have that $D=\cup_{i=1}^{s} T_{i}$ generates $I_A$ (respectively, $I_A$ up to radical). Thus, there exists a minimal system of binomial generators of $I_{A}$ (respectively, of $I_A$ up to radical) in the form $\{B({\bf u}) \mid {\bf u} \in C \subset D \}$. Let $C_i=D\cap T_{i}$ for $1 \leq i \leq s$. Thus, $C=\cup_{i=1}^{s} C_{i}$. Since $C_i\subset {\rm ker}_{\mathbb{Z}}(A_i)$, we have ${\rm span}_{\mathbb{Q}}(C_i) \subset {\rm ker}_{\mathbb{Q}}(A_i)$; therefore, ${\rm span}_{\mathbb{Q}}(C_i) \subsetneqq {\rm ker}_{\mathbb{Q}}(A)$ for every $1 \leq i \leq s$.

Conversely, suppose that there exists a minimal system of binomial generators $\{B({\bf u}_{1}),\ldots,B({\bf u}_{s}) \mid {\bf u}_{i} \in C, 1 \leq i \leq s\}$ of the toric ideal $I_{A}$ (respectively, of $I_A$ up to radical), where $C \subset {\rm ker}_{\mathbb{Z}}(A)$, and sets $C_{i}$, $1 \leq i \leq s$, such that $C=\cup_{i=1}^{s} C_{i}$ and ${\rm span}_{\mathbb{Q}}(C_i) \subsetneqq {\rm ker}_{\mathbb{Q}}(A)$ for every $1 \leq i \leq s$. The set ${\rm span}_{\mathbb{Q}}(C_i)$ is a subspace of $\mathbb{Q}^{n}$. For every $1 \leq i \leq s$, consider the matrix $M_i$ whose rows form a generating set for the orthogonal complement of the vector space ${\rm span}_{\mathbb{Q}}(C_i)$. Let $A_i$ be the set of columns of $M_i$. Since $({\rm span}_{\mathbb{Q}}(C_{i}))^{\bot \bot}={\rm ker}_{\mathbb{Q}}(A_{i})$ and also $({\rm span}_{\mathbb{Q}}(C_{i}))^{\bot \bot}={\rm span}_{\mathbb{Q}}(C_{i})$, it follows that ${\rm ker}_{\mathbb{Q}}(A_{i})={\rm span}_{\mathbb{Q}}(C_i)$ for every $1 \leq i \leq s$.
 
We claim that $I_{A}=I_{A_1}+\cdots+I_{A_s}$ (respectively, $I_{A}={\rm rad}(I_{A_1}+\cdots+I_{A_s})$) is a splitting of $I_A$ (respectively, radical splitting of $I_A$). Since $C_i\subset {\rm ker}_{\mathbb{Z}}(A)$, it follows that ${\rm ker}_{\mathbb{Z}}(A_{i}) \subset {\rm ker}_{\mathbb{Z}}(A)$. Therefore, $I_{A_i}\subset I_A$ for every $1 \leq i \leq s$, and thus ${\rm rad}(I_{A_1}+\cdots+I_{A_s}) \subset I_A$. For each generator $B({\bf u}_{i})$, where $1 \leq i \leq s$, of $I_A$ (respectively, of $I_A$ up to radical), we have that ${\bf u}_{i}\in C_j$ for $1 \leq j \leq s$. Thus $B({\bf u}_{i})\in I_{A_j}$. Therefore, every $B({\bf u}_{i})$, where $1 \leq i \leq s$, belongs to at least one of $I_{A_1}, \ldots, I_{A_s}$. We conclude that $I_A\subset I_{A_1}+\cdots+I_{A_s}$ (respectively, $I_A\subset {\rm rad}(I_{A_1}+\cdots+I_{A_s})$), and thus $I_A=I_{A_1}+\cdots+I_{A_s}$ (respectively, $I_A={\rm rad}(I_{A_1}+\cdots+I_{A_s})$). 

Note also that $I_{A_i}\not =I_A$ for $1 \leq i \leq s$. Otherwise, if $I_{A_i} =I_A$, then ${\rm ker}_{\mathbb{Z}}(A)={\rm ker}_{\mathbb{Z}}(A_{i})$, and therefore ${\rm ker}_{\mathbb{Q}}(A)={\rm span}_{\mathbb{Q}}(C_i)$, which is a contradiction. \hfill $\square$\\

The next corollary is derived directly from Theorem \ref{main2} for $s=2$.

\begin{Corollary} \label{main} The toric ideal $I_A$ is radical splittable if and only if there exists a minimal system of binomial generators $\{B({\bf u}) \mid {\bf u}\in C\subset {\rm ker}_{\mathbb{Z}}(A)\}$ of the toric ideal $I_A$ up to radical, and sets $C_{1}$ and $C_2$ such that $C=C_1 \cup C_2$, ${\rm span}_{\mathbb{Q}}(C_1)\subsetneqq {\rm ker}_{\mathbb{Q}}(A)$, and ${\rm span}_{\mathbb{Q}}(C_2)\subsetneqq {\rm ker}_{\mathbb{Q}}(A)$.   
\end{Corollary}

The following proposition establishes that $I_A$ is radical splittable whenever it is a set-theoretic complete intersection on binomials; that is, when its binomial arithmetical rank coincides with its height.

\begin{Proposition} \label{stci}
Let $K$ be a field of any characteristic and let $I_A$ be a toric ideal of height $r \geq 2$. If $I_A$ is a set-theoretic complete intersection on binomials, then $I_A$ is radical splittable. 
\end{Proposition}

{\em \noindent Proof.} Suppose that $I_A={\rm rad}(B({\bf u}_1),\ldots,B({\bf u}_r))$. Let $C=\{{\bf u}_{1}, \ldots,{\bf u}_{r}\}$, $C_{1}=\{{\bf u}_1\}$, and $C_{2}=\{{\bf u}_2,\ldots,{\bf u}_{r}\}$. Then ${\rm span}_{\mathbb{Q}}(C_{1}) \subsetneqq {\rm ker}_{\mathbb{Q}}(A)$ and ${\rm span}_{\mathbb{Q}}(C_{2}) \subsetneqq {\rm ker}_{\mathbb{Q}}(A)$, since the dimension of ${\rm ker}_{\mathbb{Q}}(A)$ is equal to $r$ and the cardinality of $C
_{i}$, $1 \leq i \leq 2$, is strictly less than $r$. By Corollary \ref{main}, the toric ideal $I_A$ is radical splittable. \hfill $\square$

\begin{Remark} {\em Note that Proposition \ref{stci} is characteristic free, but the property of $I_A$ being a set-theoretic complete intersection on binomials depends on the characteristic. For details, see \cite{BMT}.}
\end{Remark}

The next theorem shows that the corresponding statement of Proposition \ref{stci} also holds when ${\rm bar}(I_A) \leq 2{\rm ht}(I_A)-2$, provided that ${\rm ht}(I_A) \geq 3$.

\begin{Theorem} \label{BasicMinimal}
If $I_A$ is a toric ideal of height $r \geq 3$ such that its binomial arithmetical rank is less than or equal to $2r-2$, then $I_A$ is radical splittable. 
\end{Theorem}

{\em \noindent Proof.} Let $\{B({\bf u}_{1}),\ldots,B({\bf u}_{s})\}$ be a minimal system of binomial generators of $I_A$ up to radical, where $s \leq 2r-2$, and let $C=\{{\bf u}_{1},\ldots,{\bf u}_{s}\}$. Let $C_{1}=\{{\bf u}_{1},\ldots,{\bf u}_{k}\}$ and $C_{2}=\{{\bf u}_{k+1},\ldots,{\bf u}_{s}\}$, where $k=\lfloor \frac{s}{2} \rfloor$. Here, $\lfloor y \rfloor$ denotes the greatest integer less than or equal to $y$. Then, $C=C_{1} \cup C_{2}$. Since the dimension of ${\rm span}_{\mathbb{Q}}(C_{1})$ is less than or equal to $k$ and $k \leq \frac{s}{2}$, we conclude that the dimension of ${\rm span}_{\mathbb{Q}}(C_{1})$ is less than or equal to $\frac{s}{2}$. Therefore, the dimension of ${\rm span}_{\mathbb{Q}}(C_{1})$ is less than or equal to $r-1$ because $s \leq 2r-2$. But the dimension of ${\rm ker}_{\mathbb{Q}}(A)$ is equal to $r$, so the dimension of ${\rm span}_{\mathbb{Q}}(C_{1})$ is strictly smaller than the dimension of ${\rm ker}_{\mathbb{Q}}(A)$. Thus, ${\rm span}_{\mathbb{Q}}(C_{1}) \subsetneqq {\rm ker}_{\mathbb{Q}}(A)$. Note that $\frac{s}{2}<\lfloor \frac{s}{2} \rfloor+1$, and therefore, $\frac{s}{2}+1>s-\lfloor \frac{s}{2} \rfloor$. But $2r-2 \geq s$, so $r \geq \frac{s}{2}+1$ and therefore $r>s-\lfloor \frac{s}{2} \rfloor$. Since the dimension of ${\rm span}_{\mathbb{Q}}(C_{2})$ is less than or equal to $s-\lfloor \frac{s}{2} \rfloor$ and the dimension of ${\rm ker}_{\mathbb{Q}}(A)$ is equal to $r$, we have that the dimension of ${\rm span}_{\mathbb{Q}}(C_{2})$ is strictly smaller than the dimension of ${\rm ker}_{\mathbb{Q}}(A)$. Consequently, ${\rm span}_{\mathbb{Q}}(C_{2}) \subsetneqq {\rm ker}_{\mathbb{Q}}(A)$. By Corollary \ref{main}, the ideal $I_A$ is radical splittable.
\hfill $\square$

\begin{Corollary} \label{Basicbar} Let $K$ be a field of any characteristic and let $I_A$ be a toric ideal of height $r \geq 3$. If ${\rm bar}(I_A)=r+1$, then $I_A$ is radical splittable.
\end{Corollary}
{\em \noindent Proof.} Since $r \geq 3$, we have that $r+1 \leq 2r-2$. By Theorem \ref{BasicMinimal}, the ideal $I_A$ is radical splittable.
\hfill $\square$\\

Let $d_{1},\ldots,d_{n}$ be positive integers and $a_{i,j}$ be integers, where $1 \leq i \leq r$, $1 \leq j \leq n$ and for all $i=1,\ldots,r$ at least one of $a_{i,1},\ldots,a_{i,n}$ is nonzero. Let $A$ be the set of columns of the matrix \[
\left( \begin{array}{ccccccc}
 d_1 &  0 &  \cdots & 0 & a_{1,1} & \cdots & a_{r,1} \\
 0 &  d_2 &  \cdots & 0 & a_{1,2} & \cdots & a_{r,2} \\
 \cdots &  \cdots &  \cdots & \cdots & \cdots & \cdots & \cdots \\
 0 &  0 &  \cdots & d_{n} & a_{1,n} & \cdots & a_{r,n}
\end{array} \right),\]
and let $I_{A} \subset K[x_{1},\ldots,x_{n+r}]$ be the corresponding toric ideal of height $r$. The toric variety $V(I_{A}) \subset K^{n+r}$ is called {\em simplicial} if $a_{i,j}$ is nonnegative for every $1 \leq i \leq r$ and $1 \leq j \leq n$. We say that a simplicial toric variety $V(I_{A}) \subset K^{n+r}$ has a {\em full parametrization} if $a_{i,j} \neq 0$ for all $(i,j)$. 

\begin{Remark} \label{Fullbasic} {\rm Let $K$ be a field of positive characteristic. Given a simplicial toric variety $V(I_A)$ with full parametrization, we have from \cite[Theorem 1]{BMT1} that $I_A$ is a set-theoretic complete intersection on binomials, and therefore $I_A$ is radical splittable by Proposition \ref{stci}.}
\end{Remark}

The following proposition demonstrates that the toric ideal $I_A$, associated with a simplicial toric variety $V(I_A)$ that admits full parametrization, is radical splittable whenever ${\rm ht}(I_A) \geq 3$.

\begin{Proposition} \label{simplicial} Let $K$ be a field of any characteristic and let $V(I_{A}) \subset K^{n+r}$ be a simplicial toric variety with full parametrization, where ${\rm ht}(I_A)=r \geq 3$. 
Then $I_A$ is radical splittable.
\end{Proposition}

{\em \noindent Proof.} If ${\rm char}(K)>0$, then $I_A$ is radical splittable by Remark \ref{Fullbasic}. Suppose that ${\rm char}(K)=0$. If $I_A$ is a complete intersection, that is $\mu(I_A)={\rm ht}(I_A)$, then $I_A$ is splittable by \cite[Corollary 3.3]{KT3}, and therefore $I_A$ is radical splittable. Suppose that $I_A$ is not a complete intersection. Then ${\rm bar}(I_A)=r+1$ by \cite[Theorem 2]{BMT1}, and therefore $I_A$ is radical splittable by Corollary \ref{Basicbar}. \hfill $\square$

\begin{Remark} {\em Let $K$ be a field of any characteristic and let $A=\{a_{1},\ldots,a_{n}\}$ be a set of distinct positive integers, where $n \geq 4$. Then $V(I_A) \subset K^n$ is a simplicial toric variety with full parametrization and also $I_A$ has height $n-1$. By Proposition \ref{simplicial}, the toric ideal $I_A$ is radical splittable.}
\end{Remark}

Next, we provide an example of a simplicial toric variety that does not admit a full parametrization, such that $I_A$ is radical splittable.

\begin{Example} {\rm Let $A$ be the set of columns of the matrix \[
\left( \begin{array}{cccccccccc}
 3 &  0 &  0 & 2 & 1 & 0 & 0 & 2 & 1 & 1 \\
 0 &  3 &  0 & 1 & 2 & 2 & 1 & 0 & 0 & 1 \\
 0 &  0 &  3 & 0 & 0 & 1 & 2 & 1 & 2 & 1
\end{array} \right).\] Then $I_A$ is a toric ideal of height 7 which is generated up to radical by the binomials $B({\bf u}_{i})$, $1 \leq i \leq 11$, where 
\[ \begin{split}
{\bf u}_{1}=(2,1,0,-3,0,0,0,0,0,0), \\{\bf u}_{2}=(1,0,0,-2,1,0,0,0,0,0), \\{\bf u}_{3}=(1,2,0,0,-3,0,0,0,0,0),\\ 
{\bf u}_{4}=(0,2,1,0,0,-3,0,0,0,0),\\ {\bf u}_{5}=(0,1,0,0,0,-2,1,0,0,0),\\ {\bf u}_{6}=(0,1,2,0,0,0,-3,0,0,0),\\ {\bf u}_{7}=(2,0,1,0,0,0,0,-3,0,0),\\ {\bf u}_{8}=(1,0,0,0,0,0,0,-2,1,0), \\{\bf u}_{9}=(1,0,2,0,0,0,0,0,-3,0),\\ {\bf u}_{10}=(1,0,0,-1,0,0,0,-1,0,1), \\ \text{and} \ {\bf u}_{11}=(0,0,0,1,0,0,1,0,0,-2).
 \end{split} \]
From Example 3.7 in \cite{KMT}, we know that if the characteristic of the field $K$ is equal to 0, then ${\rm bar}(I_A)=11$. In this case, we may 
take $C=\{{\bf u}_1,\ldots,{\bf u}_{11}\}$, $C_{1}=\{{\bf u}_{1}, {\bf u}_{3}, {\bf u}_{4}, {\bf u}_{6}, {\bf u}_{7}, {\bf u}_{9}\}$, and $C_{2}=\{{\bf u}_{2}, {\bf u}_{5}, {\bf u}_{8}, {\bf u}_{10}, {\bf u}_{11}\}$. We have that ${\rm span}_{\mathbb{Q}}(C_1) \subsetneqq {\rm ker}_{\mathbb{Q}}(A)$ and ${\rm span}_{\mathbb{Q}}(C_2) \subsetneqq {\rm ker}_{\mathbb{Q}}(A)$, since dim$({\rm ker}_{\mathbb{Q}}(A))=7.$ From Corollary \ref{main} the ideal $I_A$ is radical splittable when the characteristic of $K$ equals 0. Next, we apply the technique described in the proof of Theorem \ref{main2} to find sets $A_1$ and $A_2$ such that ${\rm span}_{\mathbb{Q}}(C_{i})={\rm ker}_{\mathbb{Q}}(A_i)$ for $1 \leq i \leq 2$. In this way, we obtain the following two matrices: 
\[
\left( \begin{array}{cccccccccc}
 3 &  0 &  0 & 2 & 1 & 0 & 0 & 2 & 1 & 0 \\
 0 &  3 &  0 & 1 & 2 & 2 & 1 & 0 & 0 & 0 \\
 0 &  0 &  3 & 0 & 0 & 1 & 2 & 1 & 2 & 0 \\
 0 & 0 & 0 & 0 & 0 & 0 & 0 & 0 & 0 & 1
\end{array} \right),\]
\[
\left( \begin{array}{cccccccccc}
 2 &  1 &  0 & 1 & 0 & 0 & -1 & 1 & 0 & 0 \\
 -1 &  -1 &  0 & -1 & -1 & 0 & 1 & 0 & 1 & 0 \\
 0 &  -1 &  0  & 1  & 2  & 0 & 1 & 0 & 0 & 1 \\
  0 & 0 & 1 & 0 & 0 & 0 & 0 & 0 & 0 & 0 \\
   0 & 2 & 0 & 0 & 0 & 1 & 0 & 0 & 0 & 0
\end{array} \right).\]
Let $A_1$ be the set of columns of the first matrix and $A_2$ the set of columns of the second matrix. Then $I_{A}={\rm rad}(I_{A_1}+I_{A_2})$.

If $K$ is a field of characteristic $3$, then ${\rm bar}(I_A)=7$ by \cite[Example 3.7]{KMT}; namely, $I_A$ is a set-theoretic complete intersection on binomials, and therefore $I_A$ is radical splittable by Proposition \ref{stci}. In this case, $I_A$ is generated up to radical by the binomials 
$B({\bf u}_{i})$ for $ i\in \{1,3,4,6,7,9\}$, and $B({\bf v})$, where ${\bf v}=(1,1,1,0,0,0,0,0,0,-3)$. Hence, $I_{A}={\rm rad}(I_{A_1}+ \langle x_1x_2x_3-x_{10}^{3} \rangle )$. If $K$ is a field of positive characteristic $p \neq 3$, then, by \cite[Example 3.7]{KMT}, $I_A$ is generated up to radical by  the 8 binomials $B({\bf u}_{i})$, $ i\in \{1,3,4,6,7,9,10,11\}$,  and therefore $I_A$ is radical splittable by Proposition \ref{stci} and Corollary \ref{Basicbar}. 
In this case, $I_{A}={\rm rad}(I_{A_1}+ I_{A_3})$ where $I_{A_3}=\langle x_1x_{10}-x_4x_8, x_4x_7-x_{10}^2, x_1x_7-x_8x_{10} \rangle$.}
\end{Example}

\section{Applications to toric ideals of bipartite graphs}

Our aim in Section 3 is to prove that if $G$ is a bipartite graph, then ${\rm Split}_{\rm rad}(I_G)={\rm Split} (I_G)$, and that both invariants are equal to $3$ when $G$ is a complete bipartite graph. 

The radical generation of a toric ideal is a particularly difficult problem. In the literature, certain graphs or simplicial complexes (see \cite{Katsa2, matchings}) provide crucial information for identifying binomials in a minimal generating set of a toric ideal $I_A$ up to radical, or for computing the binomial arithmetical rank of $I_A$. In this section and the next, we repeatedly use techniques developed in \cite{Katsa2}.

First, we recall some material from \cite{Katsa2}. The {\em support} of a monomial ${\bf x}^{\bf u}$ of $K[x_{1},\ldots,x_{n}]$ is ${\rm supp}({\bf x}^{\bf u})=\{i \mid x_{i} \ \textrm{divides} \ {\bf x}^{\bf u}\}$. A {\em binomial} $B={\bf x}^{\bf u}-{\bf x}^{\bf v}$ is called {\em indispensable of} $I_A$ if
every system of binomial generators of $I_A$ contains $B$ or $-B$, while a {\em monomial} ${\bf x}^{\bf u}$ is called {\em indispensable of} $I_A$ if every system of binomial generators of $I_A$ contains
a binomial $B$ in which ${\bf x}^{\bf u}$ appears as one of the monomials. Let $\mathcal{T}$ be the set of all $E \subset \{1,\ldots,n\}$ for which there exists an indispensable monomial $M$ of $I_A$ such that $E ={\rm supp}(M)$. Let $\mathcal{T}_{\rm min}$ denote the set of minimal elements of $\mathcal{T}$. We associate to $I_A$ a graph $\Gamma_{A}$ with vertices the elements of $\mathcal{T}_{\rm min}$. Two distinct vertices $E, E'$ are joined by an edge if there exists a binomial $B({\bf u}) \in I_A$ such that $E={\rm supp}({\bf x}^{{\bf u}^{+}})$ and $E'={\rm supp}({\bf x}^{{\bf u}^{-}})$.

A set $M = \{T_1,\ldots,T_s\}$ is called a $\{0,1\}$-matching in $\Gamma_{A}$ if every $T_{k}, 1 \leq k \leq s$ 
is either an edge of $\Gamma_{A}$ or a singleton set of the form $\{E\}$, where $E$ is a vertex
of $\Gamma_{A}$, and also $T_{k} \cap T_{l}=\emptyset$ for all distinct indices $k,l$. Given a $\{0,1\}$-matching
$M = \{T_1,\ldots,T_s\}$ in $\Gamma_{A}$, we let ${\rm supp}(\mathcal{M})=\cup_{i=1}^{s}T_{i}$ which is a subset of the vertices
$\mathcal{T}_{\rm min}$. A $\{0,1\}$-matching $\mathcal{M}$ in $\Gamma_{A}$ is called a maximal $\{0,1\}$-matching if ${\rm supp}(\mathcal{M})$ has the maximum possible cardinality among all $\{0,1\}$-matchings. Given a maximal $\{0,1\}$-matching $\mathcal{M} = \{T_1,\ldots,T_s\}$ in $\Gamma_{A}$, we write ${\rm card}(\mathcal{M})$ for the cardinality $s$ of the set $\mathcal{M}$. By $\delta(\Gamma_{A})_{\{0,1\}}$, we denote the minimum of
the set
$$\{\textrm{card}(N) \mid N \ \textrm{is a maximal} \{0,1\}-\textrm{matching in} \ \Gamma_{A}\}.$$ By \cite[Theorem 2.6]{Katsa2}, it holds that ${\rm bar}(I_A) \geq \delta(\Gamma_{A})_{\{0,1\}}$.

The rest of this section is devoted to studying our problem in the case of toric ideals of graphs, particularly those arising from bipartite graphs.

Let $G$ be a finite, simple, connected graph in the vertex set $\{v_{1},\ldots,v_{n}\}$ with the edge set $E(G)=\{e_{1},\ldots,e_{m}\}$. For each edge $e=\{v_{i},v_{j}\}$ of $G$, we associate a vector ${\bf a}_{e} \in \{0,1\}^{n}$ defined as follows: the $i$th entry of ${\bf a}_{e}$ is $1$, the $j$th entry is $1$, and all other entries are zero. We denote the toric ideal $I_{A_{G}}$ in $K[e_{1},\ldots,e_{m}]$ by $I_G$, where $A_{G}=\{{\bf a}_{e} \mid e \in E(G)\} \subset \mathbb{N}^{n}$. 

Given an even closed walk $\gamma = (e_{i_1}, e_{i_2}, \ldots, e_{i_{2q}})$ in $G$, we define the binomial $f_{\gamma}=\prod_{k=1}^{q}e_{i_{2k-1}}-\prod_{k=1}^{q}e_{i_{2k}} \in I_{G}$. By \cite[Proposition 10.1.5]{Vil}, the ideal $I_G$ is generated by all the binomials $f_{\gamma}$, where $\gamma$ is an even closed walk of $G$. Note that if $G$ is a bipartite graph, then it contains no odd cycles, and thus $I_G$ is generated by the binomials $f_{\gamma}$, where $\gamma$ is an even cycle in $G$.

We grade $K[e_{1},\ldots,e_{m}]$ by setting ${\rm deg}_{A_{G}}(x_{i})={\bf a}_{e_i}$ for every $1 \leq i \leq m$. Given a binomial $B({\bf u}) \in I_{G}$, where ${\bf u}={\bf u}^{+}-{\bf u}^{-} \in {\rm ker}_{\mathbb{Z}}(A_{G})$, we have $${\rm deg}_{A_{G}}({\bf e}^{{\bf u}^+})={\rm deg}_{A_{G}}({\bf e}^{{\bf u}^-}).$$ For simplicity, we shall write ${\rm ker}_{\mathbb{Z}}(G)$ instead of ${\rm ker}_{\mathbb{Z}}(A_{G})$ and ${\rm ker}_{\mathbb{Q}}(G)$ instead of ${\rm ker}_{\mathbb{Q}}(A_{G})$.

In \cite[Theorem 3.2]{KT2}, it is proved that the toric ideal of $K_n$ is subgraph splittable for all $n\ge 4$, where $K_n$ denotes the complete graph with $n$ vertices. Note that for $n<4$, the ideal $I_{K_n}$ is the zero ideal. Recall that $I_{K_n}$ is subgraph splittable if there exist subgraphs $G_1$ and $G_2$ of $K_n$ such that $I_{K_n}=I_{G_1}+I_{G_2}.$ In terms of splitting numbers, this means that both ${\rm Split}(I_{K_n}), {\rm Split}_{\rm rad}(I_{K_n})$ are equal to two for $n\ge 4.$ In contrast, \cite[Theorem 4.5]{KT3}, shows that the toric ideal of the complete bipartite graph $K_{m,n}$ is not splittable.  Thus, ${\rm Split}(I_{K_{m,n}})$ is not equal to two. The aim of the next two theorems is to prove that both the splitting number and the radical splitting number of the toric ideal of a complete bipartite graph $K_{m,n}$ are equal to three, where $m,n\ge 2$ and $(m,n)\not = (2,2)$. Moreover, we prove that $I_{K_{m,n}}=I_{G_1}+I_{G_2}+I_{G_3}$, for three specific subgraphs of $K_{m,n}$.

\begin{Theorem} \label{Indispensable-Split} Let $G$ be a bipartite graph, then ${\rm Split}(I_G)={\rm Split}_{\rm rad}(I_{G})$.
\end{Theorem}

{\em \noindent Proof.} Let $V_{1}=\{x_{1},\ldots,x_{m}\}$ and $V_2=\{y_{1},\ldots,y_{n}\}$ be a bipartition of the graph $G$, such that every edge of $G$ connects a vertex in $V_1$ to a vertex in $V_2$. We denote the edge $\{x_i,y_j\}$ by $b_{ij}$. 

Suppose that ${\rm Split}_{\rm rad}(I_{G})=r$, then there exists a minimal system of binomial generators $\{B({\bf u}) \mid {\bf u}\in C\subset {\rm ker}_{\mathbb{Z}}(G)\}$ of the toric ideal $I_{G}$ up to radical, and sets $C_{i}$, $1 \leq i \leq r$ such that $C=\cup_{i=1}^{r} C_{i}$ and ${\rm span}_{\mathbb{Q}}(C_i)\subsetneqq {\rm ker}_{\mathbb{Q}}(G)$ for every $1 \leq i \leq r$. 
By \cite[Theorem 2.3]{HO}, the ideal $I_{G}$ has a unique minimal system of binomial generators $S$ consisting of all binomials of the form $f_{\gamma}$, where $\gamma$ is an even cycle of $G$ with no chord. We will show that each minimal generator $f_{\gamma}\in S$ belongs to the set $\{B({\bf u}) \mid {\bf u}\in C\subset {\rm ker}_{\mathbb{Z}}(G)\}$. Let $\gamma=(x_{i_1}, y_{j_1}, x_{i_2}, y_{j_2}, \ldots , x_{i_s}, y_{j_s})$ be the cycle  
and $f_{\gamma}={\bf e}^{{\bf v}^+}-{\bf e}^{{\bf v}^-}=b_{i_{1}j_{1}} b_{i_{2}j_{2}} \cdots  b_{i_{s}j_{s}}-
b_{i_{1}j_{s}} b_{i_{2}j_{1}} \cdots  b_{i_{s}j_{s-1}}$, where ${\bf v}={\bf v}^+-{\bf v}^- \in {\rm ker}_{\mathbb{Z}}(G)$. 
The cycle $\gamma$ has no chords; thus, the induced subgraph $W$ of $G$ on the vertex set $\{x_{i_1}, x_{i_2},  \ldots , x_{i_s}, y_{j_1}, y_{j_2}, \ldots , y_{j_s}\}$
has edge set $E(W)=\{b_{i_1j_1}, b_{i_2j_1}, b_{i_2j_2}, b_{i_3j_2}, \ldots, b_{i_sj_s}, b_{i_1j_s}\}$.
Since $f_{\gamma}$ is an indispensable binomial of $I_G$, both ${\bf e}^{{\bf v}^+}$ and ${\bf e}^{{\bf v}^-}$ are indispensable monomials; thus, we have $E={\rm supp}({\bf e}^{{\bf v}^+}) \in \mathcal{T}_{\rm min}$. By \cite[Theorem 2.4]{Katsa2}, there exists a vector ${\bf u} \in C$ such that ${\rm supp}({\bf e}^{{\bf u}^+})=E$ or ${\rm supp}({\bf e}^{{\bf u}^-})=E$.
Without loss of generality, we assume that ${\rm supp}({\bf e}^{{\bf u}^+})=E$, and also that $B({\bf u})={\bf e}^{{\bf u}^+}-{\bf e}^{{\bf u}^-}$ is irreducible. Let ${\bf e}^{{\bf u}^+}=b_{i_{1}j_{1}}^{a_{1}} b_{i_{2}j_{2}}^{a_{2}} \cdots  b_{i_{s}j_{s}}^{a_s}$. Then 
${\rm deg}_{A_{G}}({\bf e}^{{\bf u}^+})=\Sigma_{l\in [s]}a_l\epsilon_{i_l}+\Sigma_{l\in [s]}a_l\epsilon_{m+j_l}$, where $\epsilon_{t}$ is the $t$-th unit vector of $\mathbb{Z}^{m+n}$. Since ${\rm deg}_{A_{G}}({\bf e}^{{\bf u}^+})={\rm deg}_{A_{G}}({\bf e}^{{\bf u}^-})$, the variables that appear in ${\bf e}^{{\bf u}^-}$ correspond to edges of $W$. Moreover, since ${\bf e}^{{\bf u}^+}$ and ${\bf e}^{{\bf u}^-}$ have no common factor, the variables corresponding to the edges $b_{i_{1}j_{1}}, b_{i_{2}j_{2}}, \ldots, b_{i_{s}j_{s}}$ cannot appear in ${\bf e}^{{\bf u}^-}$. Thus, ${\bf e}^{{\bf u}^-}$ must be of the form 
$b_{i_{1}j_{s}}^{c_{1}} b_{i_{2}j_{1}}^{c_{2}} \cdots  b_{i_{s}j_{s-1}}^{c_s}$. Then, from ${\rm deg}_{A_{G}}({\bf e}^{{\bf u}^+})={\rm deg}_{A_{G}}({\bf e}^{{\bf u}^-})$, we have
$$\sum_{l\in [s]}a_l\epsilon_{i_l}+\sum_{l\in [s]}a_l\epsilon_{m+j_l}=\sum_{l\in [s]}c_l\epsilon_{i_l}+c_1\epsilon_{m+{j_s}}+\sum_{\substack{l \in [s] \\ l \ne 1}} c_l\epsilon_{m+j_{l-1}},$$
and therefore  $a_{1}=a_{2}=\cdots=a_{s}=c_1=c_2=\cdots =c_s$. 
Moreover, since $B({\bf u})$ is irreducible, we have $a_{1}=a_{2}=\cdots=a_{s}=c_1=c_2=\cdots =c_s=1$. Therefore, $$B({\bf u})=b_{i_{1}j_{1}} b_{i_{2}j_{2}} \cdots  b_{i_{s}j_{s}}-
b_{i_{1}j_{s}} b_{i_{2}j_{1}} \cdots  b_{i_{s}j_{s-1}}.$$ Consequently, $f_{\gamma}=B({\bf u})$. Thus, $S \subset \{B({\bf u}) \mid {\bf u}\in C\subset {\rm ker}_{\mathbb{Z}}(G)\}$. In fact $S=\{B({\bf u}) \mid {\bf u}\in C\subset {\rm ker}_{\mathbb{Z}}(G)\}$, since $\{B({\bf u}) \mid {\bf u}\in C\subset {\rm ker}_{\mathbb{Z}}(G)\}$ is a minimal system of generators of $I_{G}$ up to radical. By Theorem \ref{main2}, there exist toric ideals $I_{A_i}$ with $I_{A_i} \neq I_{G}$ for each $i$ with $1 \leq i \leq r$, such that $I_{G}=I_{A_1}+\cdots+I_{A_r}$. Consequently ${\rm Split}(I_{G}) \leq r$ and therefore ${\rm Split}(I_{G})=r$ since always ${\rm Split}_{\rm rad}(I_{G}) \leq {\rm Split} (I_{G})$. \hfill $\square$\\

A bipartite graph $G$ is called a {\em complete bipartite} graph if its vertex set can be partitioned into two subsets $V_1$ and $V_2$ such that every vertex of $V_1$ is connected to every vertex of $V_2$. It is denoted by $K_{m,n}$, where $m$ and $n$ are the numbers of vertices in $V_1$ and $V_2$ respectively. The following example provides a splitting of $I_{K_{3,3}}$.

\begin{Example} {\rm Let $K_{3,3}$ be the complete bipartite graph on the vertex set $\{x_{1},x_{2},x_{3}\} \cup \{y_{1},y_{2},y_{3}\}$ with edges $E(K_{3,3})=\{b_{ij} \mid  1 \leq i \leq 3, 1\leq j \leq 3 \}$, where $b_{ij}=\{x_{i},y_{j}\}$. Then, $I_{K_{3,3}}$ is minimally generated by the binomials $b_{21}b_{33}-b_{23}b_{31}$, $b_{22}b_{33}-b_{23}b_{32}$, $b_{21}b_{32}-b_{22}b_{31}$, $b_{11}b_{23}-b_{13}b_{21}$, $b_{12}b_{23}-b_{13}b_{22}$, $b_{11}b_{22}-b_{12}b_{21}$, $b_{11}b_{33}-b_{13}b_{31}$, $b_{12}b_{33}-b_{13}b_{32}$, and $b_{11}b_{32}-b_{12}b_{31}$. Let $G_{1}$ be the induced subgraph of $K_{3,3}$ on the vertex set $\{x_{2},x_{3}\} \cup \{y_{1},y_{2},y_{3}\}$, then $I_{G_1}$ is minimally generated by the binomials $b_{21}b_{33}-b_{23}b_{31}$, $b_{22}b_{33}-b_{23}b_{32}$, and $b_{21}b_{32}-b_{22}b_{31}$. Let $G_{2}$ be the induced subgraph of $K_{3,3}$ on the vertex set $\{x_{1},x_{2}\} \cup \{y_{1},y_{2},y_{3}\}$, then $I_{G_2}$ is minimally generated by the binomials $b_{11}b_{23}-b_{13}b_{21}$, $b_{12}b_{23}-b_{13}b_{22}$, and $b_{11}b_{22}-b_{12}b_{21}$. Let $G_{3}$ be the induced subgraph of $K_{3,3}$ on the vertex set $\{x_{1},x_{3}\} \cup \{y_{1},y_{2},y_{3}\}$, then $I_{G_3}$ is minimally generated by the binomials $b_{11}b_{33}-b_{13}b_{31}$, $b_{12}b_{33}-b_{13}b_{32}$, and $b_{11}b_{32}-b_{12}b_{31}$. Thus, $I_{K_{3,3}}=I_{G_1}+I_{G_2}+I_{G_3}$ constitutes a splitting of $I_{K_{3,3}}$.}
\end{Example}

Next, we explicitly compute ${\rm Split}_{\rm rad}(I_{K_{m,n}})$.

\begin{Theorem} \label{BipartiteSplitting} The splitting number and the radical splitting number of the toric ideal of a complete bipartite graph $K_{m,n}$ are both equal to three, where $m,n\ge 2$ and $(m,n)\not = (2,2)$.  
\end{Theorem}

{\em \noindent Proof.}  Let $V_{1}=\{x_{1},\ldots,x_{m}\}$, $V_2=\{y_{1},\ldots,y_{n}\}$ be the bipartition of the complete bipartite graph $K_{m,n}$ and $E(K_{m,n})=\{b_{ij}| 1 \leq i \leq m, 1\leq j \leq n\}$, where $b_{ij}=\{x_{i},y_{j}\}$. By \cite[Theorem 4.5]{KT3}, the ideal $I_{K_{m,n}}$ is not splittable, thus ${\rm Split}(I_{K_{m,n}}) \neq 2$, and therefore ${\rm Split}(I_{K_{m,n}}) \geq 3$. The ideal $I_{K_{m,n}}$ is minimally generated by the $2 \times 2$ minors of the matrix $M=(b_{ij})$, see \cite[Proposition 10.6.2]{Vil}. Thus, $I_{K_{m,n}}$ is minimally generated by all binomials of the form $b_{ij}b_{kl}-b_{il}b_{kj}$, such that $(b_{ij}, b_{kj}, b_{kl}, b_{il})$ is a cycle of length 4 of $K_{m,n}$.   
By the hypothesis of the theorem, at least one of the sets $V_1$, $V_2$ contains at least three vertices. Suppose $V_1$ contains three or more vertices. 
Consider the following three induced subgraphs of $K_{m,n}$: the graph $G_1$, with vertex set
$(V_1 \setminus \{x_1\})\cup V_2$; the graph $G_2$, with vertex set $(V_1 \setminus \{x_m\})\cup V_2$; and the graph $G_3$, with vertex set $\{x_1, x_m\} \cup V_2$. Let $f_{w}=b_{ij}b_{kl}-b_{il}b_{kj}$ be a minimal generator of the toric ideal of $K_{m,n}$, where $b_{ij}=\{x_{i},y_{j}\}$, $b_{kj}=\{x_{k},y_{j}\}$, $b_{kl}=\{x_{k},y_{l}\}$, and $b_{il}=\{x_{i},y_{l}\}$. Notice that $x_{i}, x_{k} \in V_{1}$ and $y_{j}, y_{l} \in V_{2}$. We distinguish the following cases. \begin{enumerate} \item If $i \neq 1$, then $f_{w}$ belongs to $I_{G_1}$. \item If $i=1$ and $k \neq m$, then $f_{w}$ belongs to $I_{G_2}$. \item If $i=1$ and $k=m$, then $f_{w}$ belongs to $I_{G_3}$.
\end{enumerate}
So, every minimal generator of $I_{K_{m,n}}$ belongs to at least one of the toric ideals $I_{G_1}$, $I_{G_2}$, or $I_{G_3}$, and therefore $I_{K_{m,n}} \subseteq I_{G_1}+I_{G_2}+I_{G_3}$. Clearly, $I_{G_i} \neq I_{K_{m,n}}$ for all $1 \leq i \leq 3$. Thus, $I_{K_{m,n}}=I_{G_1}+I_{G_2}+I_{G_3}$, and ${\rm Split}(I_{K_{m,n}})=3$. Since $K_{m,n}$ is bipartite, it follows from Theorem \ref{Indispensable-Split} that ${\rm Split}_{\rm rad}(I_{K_{m,n}})={\rm Split}(I_{K_{m,n}})$.\hfill $\square$\\

\section{The case of height 2 toric ideals}

In this section, we study the special case in which $I_A \subset K[x_{1},\ldots,x_{n}]$ is a toric ideal of height $2$. The next theorem shows that ${\rm Split}(I_A)=\mu(I_A)$ and ${\rm Split}_{\rm rad}(I_A)={\rm bar}(I_A)$ whenever ${\rm ht}(I_A)=2$.

\begin{Theorem} \label{Radicalheight2} Let $I_A$ be a toric ideal of height 2. Then the splitting number of $I_A$ (respectively, the radical splitting number of $I_A$) equals the minimal number of generators $\mu(I_A)$ of $I_A$ (respectively, the binomial arithmetical rank of $I_A$).
\end{Theorem}

{\em \noindent Proof.} Suppose that ${\rm Split}(I_{A})=r$ (respectively, ${\rm Split}_{\rm rad}(I_{A})=r$), and let $I_{A}=I_{A_1}+\cdots+I_{A_r}$ be a splitting of $I_A$ (respectively, let $I_{A}={\rm rad}(I_{A_1}+\cdots+I_{A_r})$ be a radical splitting of $I_A$). Since $I_{A_i} \subsetneqq I_A$, we have that ${\rm ht}(I_{A_i})<{\rm ht}(I_A)$, and therefore ${\rm ht}(I_{A_i})=1$ for every $1 \leq i \leq r$. Thus, each $I_{A_i}$ is a principal ideal, and let $T_{i}$ be a minimal binomial generating set of $I_{A_i}$, which is a singleton. Then $D=\cup_{i=1}^{r}T_{i}$ generates $I_A$ (respectively, $D$ generates $I_A$ up to radical). So $r \geq \mu(I_A)$ (respectively, $r \geq {\rm bar}(I_A)$), and therefore $\mu(I_A)=r$ (respectively, ${\rm bar}(I_A)=r$), since ${\rm Split}(I_{A}) \leq \mu(I_A)$ (respectively, ${\rm Split}_{\rm rad}(I_{A}) \leq {\rm bar}(I_A)$). \hfill $\square$\\

The next corollary provides a necessary and sufficient condition for $I_A$ to be radical splittable.

\begin{Corollary} \label{height2}
Let $K$ be a field of any characteristic and let $I_A$ be a toric ideal of height 2. Then $I_A$ is radical splittable if and only if $I_A$ is a set-theoretic complete intersection on binomials. 
\end{Corollary}

{\em \noindent Proof.} Since ${\rm Split}_{\rm rad}(I_{A})={\rm bar}(I_A)$ by Theorem \ref{Radicalheight2}, we have that $I_A$ is radical splittable if and only if ${\rm bar}(I_{A})=2$; in other words $I_{A}$ is a set-theoretic complete intersection on binomials. \hfill $\square$

\begin{Remark} \label{Non-Stci} {\rm Suppose that ${\rm char}(K)=0$ and that the height of $I_A$ is 2. If $I_A$ is not a complete intersection, then $I_A$ is not a set-theoretic complete intersection on binomials by \cite[Theorem 4]{BMT}, and therefore ${\rm Split}_{\rm rad}(I_A) \geq 3$.}
\end{Remark}

The next corollary provides bounds for ${\rm Split}_{\rm rad}(I_A)$ whenever $V(I_A)$ is a simplicial toric variety.

\begin{Corollary} \label{Simplicial-2}
Let $K$ be a field of any characteristic, and let $V(I_A)$ be a simplicial toric variety, where $I_A$ is a toric ideal of height 2. Then $2 \leq {\rm Split}_{\rm rad}(I_A) \leq 3$. 
\end{Corollary}

{\em \noindent Proof.} By \cite[Theorem 3]{BMT1}, we have $2 \leq {\rm bar}(I_A) \leq 3$, and therefore $2 \leq {\rm Split}_{\rm rad}(I_A) \leq 3$. \hfill $\square$\\

Next, we provide an example of a simplicial toric variety $V(I_A)$ with full parametrization, such that ${\rm ht}(I_A)=2$ and ${\rm Split}_{\rm rad}(I_A)=3$ when the characteristic of $K$ is zero.

\begin{Example} \label{Simplicial-Example} {\rm Let $A$ be the set of columns of the matrix \[
\left( \begin{array}{cccc}
 2 &  1 &  2 & 0   \\
 3 &  0 &  2 & 5 
\end{array} \right).\] Then $I_A$ is a toric ideal of height 2 which is minimally generated by the binomials $x_1^3-x_{2}^2x_{3}^2x_{4}$, $ x_3^3-x_1^2x_2^2$, and $x_1x_3-x_2^4x_4$. By Corollary \ref{Simplicial-2}, it holds $2 \leq {\rm Split}_{\rm rad}(I_A) \leq 3$. If ${\rm char}(K)>0$, then $I_A$ is radical splittable by Remark \ref{Fullbasic}, and therefore ${\rm Split}_{\rm rad}(I_A)=2$. If ${\rm char}(K)=0$, then ${\rm Split}_{\rm rad}(I_A) \geq 3$ by Remark \ref{Non-Stci}, and therefore ${\rm Split}_{\rm rad}(I_A)=3$.
}
\end{Example}

We now provide an example of a non-simplicial toric variety $V(I_A)$ such that ${\rm ht}(I_A)=2$ and ${\rm Split}_{\rm rad}(I_A)=3$, in any characteristic.

\begin{Example} {\rm Let $A$ be the set of columns of the matrix \[
\left( \begin{array}{ccccc}
 33 &  -4 &  1 & 0 & 1  \\
 23 &  -3 &  0 & -9 & 0 \\
 -34 &  6 &  0 & 3 & 1
\end{array} \right).\] Then $I_A$ is a toric ideal of height 2 which is minimally generated by the binomials $x_5^{21}-x_1^3x_2^{20}x_{3}^2x_4$, $x_1^3x_2^{17}x_4^2-x_3^{25}x_5^6$, $x_3^{52}- x_1^3x_2^{14}x_4^3x_5^9$, and $x_2^3x_3^{27}- x_4x_5^{15}$. We have that $I_A={\rm rad}(x_5^{21}-x_1^3x_2^{20}x_{3}^2x_4, x_3^{52}- x_1^3x_2^{14}x_4^3x_5^9,x_1^3x_2^{17}x_4^2-x_3^{25}x_5^6)$, since $(x_2^3x_3^{27}- x_4x_5^{15})^2$ belongs to the ideal generated by the binomials $x_5^{21}-x_1^3x_2^{20}x_{3}^2x_4$, $x_1^3x_2^{17}x_4^2-x_3^{25}x_5^6$, and $x_3^{52}- x_1^3x_2^{14}x_4^3x_5^9$. Thus ${\rm bar}(I_A) \leq 3$. Then $\mathcal{T}_{{\rm min}}=\{\{3\},\{5\}, \{1,2,4\}\}$ and $\Gamma_{A}$ consists only of the three vertices of $\mathcal{T}_{{\rm min}}$. Thus $\delta(\Gamma_{A})_{\{0,1\}}=3$ and therefore ${\rm bar}(I_A) \geq 3$. Consequently ${\rm bar}(I_A)=3$, so the radical splitting number of $I_A$ equals 3 in any characteristic.}
\end{Example}

Next, we provide two classes of toric ideals for which the splitting number can be arbitrarily large, while the radical splitting number remains constant. The first class determines a simplicial toric variety, while the second determines a non-simplicial toric variety.

Let $a<b<d$ be positive integers with ${\rm gcd}(a,b,d)=1$ and $a+b=d$. A {\em symmetric monomial curve} in $\mathbb{P}^{3}$ is given parametrically by $$(u^{d},u^{d-a}v^{a},u^{d-b}v^{b},v^{d}).$$ Let $S$ be the set of columns of the $2 \times 4$ matrix \[
\left( \begin{array}{cccc}
 d &  d-a &  d-b & 0    \\
 0 &  a &  b & d 
\end{array} \right). \] Then $I_S$ is minimally generated (see \cite{BSV}) by the $b-a+2$ binomials $x_{1}x_{4}-x_{2}x_{3}$ and $x_{1}^{b-a-i}x_{3}^{a+i}-x_{2}^{b-i}x_{4}^{i}$, for every $0 \leq i \leq b-a$. By Theorem \ref{Radicalheight2}, the splitting number of $I_S$ equals $b-a+2$.

Note that $I_S=I_{A}$ where $A$ is the set of columns of the $2 \times 4$ matrix \[
\left( \begin{array}{cccc}
 1 &  1 &  1 & 1    \\
 0 &  a &  b & a+b 
\end{array} \right). \] 

\begin{Remark} {\rm Let $I_A$ be the toric ideal of a symmetric monomial curve in $\mathbb{P}^{3}$. Since $V(I_S)$ is a simplicial toric variety with full parametrization, it follows from Remark \ref{Fullbasic} that the radical splitting number of $I_{A}$ is $2$ whenever ${\rm char}(K)>0$. If ${\rm char}(K)=0$, then ${\rm Split}_{\rm rad}(I_A) \geq 3$ by Remark \ref{Non-Stci}, and hence ${\rm Split}_{\rm rad}(I_A)=3$.}
\end{Remark}

The Lawrence lifting $\Lambda(A)$ of $A$ is the set of columns of the $6 \times 8$ matrix \[
\left( \begin{array}{cccccccc}
 1 &  1 &  1 & 1 & 0 & 0 & 0 & 0    \\
 0 &  a &  b & a+b & 0 & 0 & 0 & 0 \\
 1 & 0 & 0 & 0 & 1 & 0 & 0 & 0 \\
 0 & 1 & 0 & 0 & 0 & 1 & 0 & 0 \\
 0 & 0 & 1 & 0 & 0 & 0 & 1 & 0 \\
 0 & 0 & 0 & 1 & 0 & 0 & 0 & 1 
\end{array} \right). \]
The toric ideal $I_{\Lambda(A)} \subset K[x_{1},\ldots,x_{4},y_{1},\ldots,y_{4}]$ is generated by all binomials ${\bf x}^{{\bf u}^{+}}{\bf y}^{{\bf u}^{-}}-{\bf x}^{{\bf u}^{-}}{\bf y}^{{\bf u}^{+}}$ where ${\bf u}={\bf u}^{+}-{\bf u}^{-} \in {\rm ker}_{\mathbb{Z}}(A)$. Our aim is to compute ${\rm Split}_{\rm rad}(I_{\Lambda(A)})$.

Given a toric ideal $I_A \subset K[x_{1},\ldots,x_{n}]$, a binomial ${\bf x}^{{\bf u}^+}-{\bf x}^{{\bf u}^-} \in I_A$ is called {\em primitive} if there is no other binomial ${\bf x}^{{\bf v}^+}-{\bf x}^{{\bf v}^-} \in I_A$ such that ${\bf x}^{{\bf v}^+}$ divides ${\bf x}^{{\bf u}^+}$ and ${\bf x}^{{\bf v}^-}$ divides ${\bf x}^{{\bf u}^-}$. The set of the primitive binomials
is finite, it is called the Graver basis of $I_A$ and is denoted by ${\rm Gr}(I_A)$. An irreducible binomial ${\bf x}^{{\bf u}^+}-{\bf x}^{{\bf u}^-} \in I_A$ is called a {\em circuit} of $I_A$ if there exists no binomial ${\bf x}^{{\bf v}^+}-{\bf x}^{{\bf v}^-} \in I_A$ such that ${\rm supp}({\bf x}^{{\bf v}^+}) \cup {\rm supp}({\bf x}^{{\bf v}^-}) \subsetneqq {\rm supp}({\bf x}^{{\bf u}^+}) \cup {\rm supp}({\bf x}^{{\bf u}^-})$. Let $\mathcal{C}$ denote the collection of subsets $E \subset \{1,\ldots,n\}$ such that ${\rm supp}({\bf x}^{{\bf u}^{+}})=E$ or ${\rm supp}({\bf x}^{{\bf u}^-})=E$, where $B({\bf u})$ is a circuit of $I_A$. Let $\mathcal{C}_{\rm min}$ denote the set of the minimal elements of $\mathcal{C}$. By \cite[Proposition 2.2]{Katsa2}, the set $\mathcal{T}_{\rm min}$ is equal to $\mathcal{C}_{\rm min}$.

\begin{Proposition} \label{Lawrence}
For the Lawrence lifting $\Lambda(A)$ of $A$, the following statements hold: \begin{enumerate} \item The splitting number of $I_{\Lambda(A)}$ is $a+b+2$. \item The radical splitting number of $I_{\Lambda(A)}$ is $4$.
\end{enumerate}
\end{Proposition}

{\em \noindent Proof.} Let ${\bf v}_{i}=(b-i,i-a-b,i,a-i)$ for every $0 \leq i \leq a+b$, and let ${\bf w}=(1,-1,-1,1)$. By \cite[Lemma 2.2]{Nai}, the Graver basis of $I_A$ is $${\rm Gr}(I_A)=\{B({\bf v}_{i}) \mid  0 \leq i \leq a+b\} \cup \{B({\bf w})\}.$$ Then $I_{\Lambda(A)}$ is minimally generated by the $a+b+2$ binomials $B({\bf w},-{\bf w})={\bf x}^{{\bf w}^{+}}{\bf y}^{{\bf w}^{-}}-{\bf x}^{{\bf w}^{-}}{\bf y}^{{\bf w}^{+}}$ and $B({\bf v}_{i},-{\bf v}_{i})={\bf x}^{{\bf v}_{i}^{+}}{\bf y}^{{\bf v}_{i}^{-}}-{\bf x}^{{\bf v}_{i}^{-}}{\bf y}^{{\bf v}_{i}^{+}}$ for every $0 \leq i \leq a+b$.

(1) By Theorem \ref{Radicalheight2}, the splitting number of $I_{\Lambda(A)}$ is $a+b+2$.\\ 

(2) Since $${\rm Gr}(I_A)=\{B({\bf v}_{i}) \mid  0 \leq i \leq a+b\} \cup \{B({\bf w})\},$$ we have that $I_A$ has 4 circuits, namely $x_{2}^{b}-x_{1}^{b-a}x_{3}^{a}$, $x_{2}^{a+b}-x_{1}^{b}x_{4}^{a}$, $x_{3}^{b}-x_{2}^{a}x_{4}^{b-a}$, and $x_{3}^{a+b}-x_{1}^{a}x_{4}^{b}$. So $I_{\Lambda(A)}$ has 4 circuits, namely $x_{2}^{b}y_{1}^{b-a}y_{3}^{a}-x_{1}^{b-a}x_{3}^{a}y_{2}^{b}$, $x_{2}^{a+b}y_{1}^{b}y_{4}^{a}-x_{1}^{b}x_{4}^{a}y_{2}^{a+b}$, $x_{3}^{b}y_{2}^{a}y_{4}^{b-a}-x_{2}^{a}x_{4}^{b-a}y_{3}^{b}$, and $x_{3}^{a+b}y_{1}^{a}y_{4}^{b}-x_{1}^{a}x_{4}^{b}y_{3}^{a+b}$. By \cite[Proposition 8.7]{ES}, $I_{\Lambda(A)}$ is generated up to radical by its circuits. Thus ${\rm bar}(I_{\Lambda(A)}) \leq 4$. Let $E_{1}={\rm supp}(x_{2}^{b}y_{1}^{b-a}y_{3}^{a})$, $E_{2}={\rm supp}(x_{2}^{a+b}y_{1}^{b}y_{4}^{a})$, $E_{3}={\rm supp}(x_{3}^{b}y_{2}^{a}y_{4}^{b-a})$, $E_{4}={\rm supp}(x_{3}^{a+b}y_{1}^{a}y_{4}^{b})$, $E_{5}={\rm supp}(x_{1}^{b-a}x_{3}^{a}y_{2}^{b})$, $E_{6}={\rm supp}(x_{1}^{b}x_{4}^{a}y_{2}^{a+b})$, $E_{7}={\rm supp}(x_{2}^{a}x_{4}^{b-a}y_{3}^{b})$, and $E_{8}={\rm supp}(x_{1}^{a}x_{4}^{b}y_{3}^{a+b})$. Clearly, $\mathcal{C}_{\rm min}=\{E_{i} \mid 1 \leq i \leq 8\}$. Then the graph $\Gamma_{\Lambda(A)}$ consists of four edges, namely $\{E_{1},E_{2}\}$, $\{E_{3},E_{4}\}$, $\{E_{5},E_{6}\}$, and $\{E_{7},E_{8}\}$. Thus, $\delta(\Gamma_{\Lambda(A)})_{\{0,1\}}=4$, and therefore ${\rm bar}(I_{\Lambda(A)}) \geq 4$. Consequently, ${\rm bar}(I_{\Lambda(A)})=4$, and thus ${\rm Split}_{\rm rad}(I_{\Lambda(A)})=4$. \hfill $\square$\\

We conclude the paper by presenting a class of toric ideals $I_A$ of height $2$ for which both the splitting number and the radical splitting number attain arbitrarily large values. We emphasize that cyclic configurations, known for their extremal properties (see, e.g., \cite{MM}), play a key role in our construction.

Given an integer $t$, define  \[ {\bf a}_t=
\left( \begin{array}{c}
 1\\ t \\ t^2\\ \vdots \\t^{2d-2}
\end{array} \right) \in \mathbb{Z}^{2d-1},\] where $d \geq 2$ is an integer.
Let $A$ denote the cyclic configuration formed by the columns of the $(2d-1) \times (2d+1)$ Vandermonde matrix \[ 
\left( \begin{array}{cccc}
 {\bf a}_{t_1} & {\bf a}_{t_2} & \hdots & {\bf a}_{t_{2d+1}} 
\end{array} \right), \] where $t_1< t_2< \cdots <t_{2d+1}$ are integers. Any subset of $A$ consisting of $2d-1$ vectors is linearly independent, so the dimension of ${\rm ker}_{\mathbb{Q}}(A)$ equals $2$, and therefore ${\rm ht}(I_A)=2$. 

\begin{Proposition} \label{cyclic} Let $I_A$ be the toric ideal of a cyclic configuration
of height two.
The splitting number and the radical splitting number of $I_{A}$ can attain arbitrarily large values.
\end{Proposition}

{\em \noindent Proof.} Any subset of $A$ consisting of $2d$ vectors defines a circuit of $I_A$. Consequently, there are exactly $2d+1$ circuits in $I_A$. By \cite[Proposition 8.7]{ES}, $I_{A}$ is generated up to radical by its circuits, thus ${\rm bar}(I_{A}) \leq 2d+1$. Consider an increasing sequence $t_{i_1}< t_{i_2}< \cdots <t_{i_{2d}}$ obtained by omitting a single element from $\{t_{1},\ldots,t_{2d+1}\}$. Then, by Cramer's rule, we obtain $$\sum_{j=1}^{2d} (-1)^j \det ({\bf a}_{t_{i_1}},{\bf a}_{t_{i_2}}, \ldots ,{\bf a}_{t_{i_{j-1}}}, {\bf a}_{t_{i_{j+1}}}, \ldots, {\bf a}_{t_{i_{2d}}}) \epsilon_{i_j}$$ is a scalar multiple of a circuit, where $\epsilon_i$ denotes the $i$-th unit vector of $\mathbb{Z}^{2d+1}$. Recall that the determinant $\det ({\bf a}_{t_{i_1}},{\bf a}_{t_{i_2}}, \ldots ,{\bf a}_{t_{i_{j-1}}}, {\bf a}_{t_{i_{j+1}}}, \ldots, {\bf a}_{t_{i_{2d}}}) $  is given by $$\prod_{1 \leq k<l \leq 2d}^{k,l\not=j}(t_{i_l}-t_{i_k})>0.$$
Thus the signs of the nonzero components of a circuit alternate. Then every element of $\mathcal{C}_{\rm min}$ has exactly $d$ elements and can be represented as $\{k_1, k_2, \ldots, k_d\}$,
where $1\leq k_1< k_2< \cdots< k_d\leq 2d+1$ and $k_{j+1}-k_j=2$ for all but at most one index
$j$, for which $k_{j+1}-k_j=3$. Each circuit contributes two elements in $\mathcal{C}_{\rm min}$, and each element of $\mathcal{C}_{\rm min}$ appears in exactly two circuits. Therefore, the graph $\Gamma_A$ is a cycle with $2d+1$ vertices and $2d+1$ edges. Thus, $\delta(\Gamma_{A})_{\{0,1\}}=d+1$, and therefore ${\rm bar}(I_{A}) \geq d+1$. Consequently,  $d+1 \leq {\rm Split}_{\rm rad}(I_{A}) \leq 2d+1$, since ${\rm Split}_{\rm rad}(I_{A})={\rm bar}(I_A)$ by Theorem \ref{Radicalheight2}. Since ${\rm Split}_{\rm rad}(I_{A}) \leq {\rm Split}(I_{A})$, we conclude that $d+1 \leq {\rm Split}(I_{A})$. \hfill $\square$\\

\noindent \textbf{Funding} The authors did not receive support from any organization for the submitted work.\\

\noindent \textbf{Data availability} Data sharing not applicable to this article as no datasets were generated or analysed during the current study.\\

\noindent \textbf{Declaration}\\

\noindent \textbf{Conflict of interest}  The authors declare that there are no conflict of interest that are directly or indirectly related to the work submitted for publication.

\end{document}